\documentclass[epj,nopacs]{svjourArXiv}

\usepackage{amsmath,amssymb}

\usepackage{stackrel}

\newcommand{\Q}{\mathbb{Q}}

\newcommand{\rationals}{\Q}

\newcommand{\reals}{\mathbb{R}}

\makeatletter
\newif\if@borderstar
\def\bordermatrix{\@ifnextchar*{%
\@borderstartrue\@bordermatrix@i}{\@borderstarfalse\@bordermatrix@i*}%
}
\def\@bordermatrix@i*{\@ifnextchar[{\@bordermatrix@ii}{\@bordermatrix@ii[() ]}}
\def\@bordermatrix@ii[#1]#2{%
\begingroup
\m@th\@tempdima8.75\p@\setbox\z@\vbox{%
\def\cr{\crcr\noalign{\kern 2\p@\global\let\cr\endline }}%
\ialign {$##$\hfil\kern 2\p@\kern\@tempdima & \thinspace %
\hfil $##$\hfil && \quad\hfil $##$\hfil\crcr\omit\strut %
\hfil\crcr\noalign{\kern -\baselineskip}#2\crcr\omit %
\strut\cr}}%
\setbox\tw@\vbox{\unvcopy\z@\global\setbox\@ne\lastbox}%
\setbox\tw@\hbox{\unhbox\@ne\unskip\global\setbox\@ne\lastbox}%
\setbox\tw@\hbox{%
$\kern\wd\@ne\kern -\@tempdima\left\@firstoftwo#1%
\if@borderstar\kern2pt\else\kern -\wd\@ne\fi%
\global\setbox\@ne\vbox{\box\@ne\if@borderstar\else\kern 2\p@\fi}%
\vcenter{\if@borderstar\else\kern -\ht\@ne\fi%
\unvbox\z@\kern-\if@borderstar2\fi\baselineskip}%
\if@borderstar\kern-2\@tempdima\kern2\p@\else\,\fi\right\@secondoftwo#1 $%
}\null \;\vbox{\kern\ht\@ne\box\tw@}%
\endgroup}
\makeatother 

\newcommand{\absoluteval}[1]{\left\vert#1\right\vert}

\newcommand{\covat}[3]{\Cov_{#1}\left(#2,#3\right)}

\newcommand{\expectat}[2]{{\E}_{#1}\left[#2\right]}

\newcommand{\expof}[1]{\exp\left(#1\right)}

\newcommand{\kcoshof}[1]{\Kcosh\left(#1\right)}
\newcommand{\kexpof}[1]{\Kexp\left(#1\right)}
\newcommand{\klnof}[1]{\Kln\left(#1\right)}
\newcommand{\kplus}{\stackbin[]{\kappa}{\oplus}}
\newcommand{\ktimes}{\stackbin[\kappa]{}{\otimes}}
\newcommand{\lnof}[1]{\ln\left(#1\right)}
\newcommand{\logof}[1]{\log\left(#1\right)}

\newcommand{\setof}[2]{\left\{#1 : #2 \right\}}
\newcommand{\set}[1]{\left\{#1\right\}}

\newcommand{\spanof}[1]{\Span\left(#1\right)}

\newcommand{\euler}{\mathrm{e}}

\DeclareMathOperator{\Cov}{Cov}

\DeclareMathOperator{\E}{E}

\DeclareMathOperator{\Kcosh}{\cosh_{\kappa}}
\DeclareMathOperator{\Kexp}{\exp_{\kappa}}
\DeclareMathOperator{\Kln}{\ln_{\kappa}}

\DeclareMathOperator{\Span}{Span}

\def\cocoa{{\hbox{\rm C\kern-.13em o\kern-.07em C\kern-.13em o\kern-.15em A}}}

\newcommand{\densities}{\mathcal M_{\ge}}
\newcommand{\pdensities}{\mathcal M_{>}}
\newcommand{\sdensities}{\mathcal M^1}
\begin{document}
\title{$\kappa$-exponential models from the geometrical viewpoint}
\author{Giovanni Pistone}
\institute{Politecnico di Torino \email{giovanni.pistone@polito.it}}
\date{\today}
%\thanks{}
%
%
\abstract{We discuss the use of Kaniadakis' $\kappa$-exponential in the construction of a statistical manifold modelled on Lebesgue spaces of real random variables. Some algebraic features of the deformed exponential models are considered. A chart is defined for each strictly positive densities; every other strictly positive density in a suitable neighborhood of the reference probability is represented by the centered $\Kln$ likelihood.}
\maketitle
\section{Introduction}
G. Kaniadakis \cite{kaniadakis:2001PhA,kaniadakis:2001PhLA,kaniadakis:2002PhRE,kaniadakis:2005PhRE}, based on arguments from Statistical Physics and Special Relativity, has defined the \emph{$\kappa$-deformed exponential} for each $x \in \reals$ and $-1 < \kappa < 1$ to be
\begin{equation}
  \kexpof x = \expof{\int_0^x \frac{dt}{\sqrt{1+\kappa^2 t^2}}},
\end{equation}
with special cases
\begin{align}
\kexpof x &= \begin{cases}
\left(\kappa x + \sqrt{1 + \kappa^2 x^2}\right)^{\frac1\kappa}, &\text{if $\kappa \ne 0$,}\\
\exp x, &\text{if $\kappa = 0$,}
\end{cases}
\end{align}
and derivation formul\ae
\begin{align}
  \Kexp'(x) &= (1+\kappa^2 x^2)^{-1/2} \kexpof x > 0 \\
  \Kexp''(x) &= \frac{\sqrt{1+\kappa^2x^2} - \kappa^2 x}{1+\kappa^2x^2} \Kexp'(x)>0
\end{align}
For each $\kappa\ne0$, $y = (\kexpof x)^\kappa$ and $x$ are related by the polynomial equation
\begin{equation} \label{eq:master}
  y^2 - 2\kappa xy - 1 = 0
\end{equation}
Therefore, the graph of $(\Kexp)^\kappa$ is the upper branch of a hyperbola:
\begin{equation}
  \label{eq:hyperbola}
  x = \frac1{2\kappa}\left(y - \frac1y\right), \quad y>0.
\end{equation}

For each given $\kappa$, the function $\Kexp$ maps $\reals$ onto $\reals_>$, it is strictly increasing and it is strictly convex. If $\kappa \ne 0$, its inverse function is
\begin{equation}\label{eq:klog}
  \klnof y = \frac{y^\kappa - y^{-\kappa}}{2\kappa}, \quad y > 0,
\end{equation}
with derivative
\begin{equation}
  \label{eq:Klnderiv}
  \Kln'(y) = \frac{y^\kappa+y^{-\kappa}}{2}\frac1y,
\end{equation}
and it is called \emph{$\kappa$-deformed logarithm}. The function $\Kln$ maps $\reals_>$ unto $\reals$, is strictly increasing and is strictly concave. Both deformed exponential and logarithm functions $\Kexp$ and $\Kln$ reduce to the ordinary $\exp$ and $\ln$ functions when $\kappa \to 0$. Moreover, the algebraic properties of the exponential and logarithmic function are partially preserved, because
\begin{equation} \label{eq:kexpofminus}
  \kexpof x \kexpof{-x} = 1, \quad \klnof y + \klnof {y^{-1}} = 0.
\end{equation}

Kaniadakis \cite{kaniadakis:2002PhRE} defines, among others, two commutative group operations $(\reals, \kplus)$ and $(\reals_>,\ktimes)$ in such a way that $\Kexp$ is a group isomorphism from $(\reals,+)$ onto $(\reals_>,\ktimes)$ and also from $(\reals,\kplus)$ onto $(\reals_>,\times)$:
\begin{align}
  \kexpof{x_1 \kplus x_2} &= \kexpof{x_1} \kexpof{x_2}, \\
  \kexpof{x_1 + x_2} &= \kexpof{x_1} \ktimes \kexpof{x_2},
\end{align}
and, equivalently,
\begin{align}
  \klnof{y_1 \ktimes y_2} &= \klnof{y_1} + \klnof{y_2}, \\
  \klnof{y_1 y_2} &= \klnof{y_1} \kplus \klnof{y_2}.
\end{align}
Therefore, the binary operations $\kplus$ and $\ktimes$ are defined by
\begin{align}
  x_1 \kplus x_2 &= \klnof{\kexpof{x_1}\kexpof{x_2}} \label{eq:kplus-def}\\
  y_1 \ktimes y_2 &=\kexpof{\klnof{y_1} + \klnof{y_2}} \label{eq:kprod-def}
\end{align}
Because of \eqref{eq:kexpofminus}, the deformed operation have the same inverse that the usual operations:
\begin{equation}
  \label{eq:kinverse}
  x_1 \kplus (-x_2) = 0, \quad y_1 \ktimes y_2^{-1} = 1.
\end{equation}
The operation $\ktimes$ is defined on positive real numbers. However, \eqref{eq:kprod-def} can be extended by continuity to non-negative real numbers:
\begin{equation}
  \label{eq:kprod-by-zero}
  0 \ktimes y_2= y_1 \ktimes 0 = 0 \ktimes 0 = 0
\end{equation}

We want to derive defining relations for the $\kappa$-deformed operations in polynomial form. This is obtained by repeated use of \eqref{eq:master} followed by algebraic elimination of the unwanted indeterminate. Symbolic computations have been done with \cite{CocoaSystem}. First, we want to find $x=x_1\kplus x_2$, i.e. such that $\kexpof x = \kexpof{x_1}\kexpof{x_2}$.
Let $y_1 = (\kexpof{x_1})^\kappa$ and $y_2 = (\kexpof{x_2})^\kappa$. As
\begin{equation}
  (\kexpof x)^\kappa = (\kexpof{x_1}\kexpof{x_2})^\kappa = y_1y_2,
\end{equation}
Equation \eqref{eq:master} gives a system of three quadratic equations in the indeterminates $x_1,x_2,y_1,y_2,x,\kappa$. Algebraic elimination of $y_1,y_2$ gives the polynomial equation
\begin{equation}
x^{4} - 2\left(2 \kappa^{2}x_{1}^{2}x_{2}^{2} + x_{1}^{2} + x_{2}^{2}\right)x^{2} + \left(x_{1}^{2} - x_{2}^{2}\right)^2 = 0, 
\end{equation}
whose solution gives
\begin{equation}
   x_1 \kplus x_2 = x_1 \sqrt{1+\kappa^2 x_2^2} + x_2 \sqrt{1+\kappa^2 x_1^2}.
\end{equation}
We will use later on the derivation formula
\begin{equation}\label{eq:11}
  \frac{\partial (x_1 \kplus x_2)}{\partial x_1} = \sqrt{1+\kappa^2 x_2^2} - \frac{\kappa^2 x_1x_2}{\sqrt{1+\kappa^2 x_1^2}}.
\end{equation}

Second, we want to find $z = \left(y_1\ktimes y_2\right)^\kappa$. As before, $y_1 = (\kexpof{x_1})^\kappa$, $y_2 = (\kexpof{x_2})^\kappa$, while
\begin{equation}
  z = (\kexpof{x_1 + x_2})^\kappa.
\end{equation}
Equation \eqref{eq:master} gives three quadratic equations in the indeterminates $x_1,x_2,y_1,y_2,z,\kappa$. Elimination of $x_1,x_2$ gives the polynomial equation
\begin{equation}
  y_{1}y_{2}z^{2} +(1 - y_{1}y_{2})(y_{1} + y_{2}) z - y_{1}y_{2} = 0
\end{equation}
It is remarkable that this equation does not depend on $\kappa$. An explicit solution is obtained by solving the quadratic equation. A possibly more suggestive solution is obtained as follows. First, we reduce to the monic equation
\begin{equation}
  \label{eq:kprod-2}
  z^{2} + \left(1 - \frac1{y_{1}y_{2}}\right)(y_{1} + y_{2}) z - 1 = 0
\end{equation}
and denote the two solutions as $z >0$ and $-1/z$. Therefore,
\begin{equation} \label{eq:15}
z - \frac1z = \left(y_1 - \frac1{y_1}\right) + \left(y_2 - \frac1{y_2}\right).
\end{equation}

The $\kappa$-logarithm defined in \eqref{eq:klog} is reminiscent of a family of transformation well known in Applied Statistics under the name of Box-Cox transformation \cite{box|cox:64} or power transform. For data vector $y_1,\dots, y_n$, with $y_i > 0$ for all $i$, the power transform is:
\begin{equation}
  y_i^{(\lambda)} \propto \frac{y_i^\lambda -1}{\lambda}
\end{equation}
The parameter $\lambda$ is to be estimated in order to get the best fit to the Normal distribution of the transformed data vector. When compared with the power transform, the $\kappa$-deformed logarithm $x = \klnof y$ has the extra feature of the symmetry induced by the term $-y^{-\kappa}$ and it would be interesting to study it from the point of view of transformations to normality. We are not further discussing this issue here.  

J. Naudts \cite{naudts:2002PhysA,naudts:2004-332,naudts:2004JIPAM,naudts:2008} has presented a general discussion of a class of deformed exponentials that contains $\Kln$ and $\Kexp$, with applications in Information Theory and Statistical Physics, based on a notion of generalised entropy. We are not discussing generalised entropies in this paper. The purpose of the present paper is to extend to $\kappa$-deformations the non-parametric and geometric approach to statistical manifolds as it was developed by the Author and co-workers in \cite{pistone|sempi:95,gibilisco|pistone:98,pistone|rogantin:99,cena|pistone:2007,pistone:2008}. Such an approach is designed to present in a non-parametric and fully geometric way work by S-i. Amari on Information Geometry, see \cite{amari:82} and the joint monograph with H. Nagaoka \cite{amari|nagaoka:2000}. According to those references, the the power transform is applied to probability densities and re-named $\alpha$-embedding, $\alpha=1-2/\lambda$. 

To explain the idea of $\alpha$-embedding, we consider a general setting. Let $\Omega$ be any set, $\mathcal F$ a $\sigma$-algebra of subsets, $\mu$ a reference probability measure, e.g. the uniform distribution. A density $p$ of the measure space $(\Omega,\mathcal F, \mu)$ is a non-negative random variable such that $\int p \,d\mu = 1$. We write the expected value of a random variable $U$ with respect to the probability measure $p \cdot \mu$ as $\expectat p U = \int Up\,d\mu$. For $a > 1$, the mapping $p \mapsto p^{1/a}$ maps $p$ into the unit ball of the Lebesgue space $L^a(\Omega,\mathcal F, \mu)$. Other parametrization of the same setting are: $\lambda = 1/a$ (Box-Cox), $\kappa = 1/a$ (Kaniadakis), $\alpha=(a-2)/a$ (Amari). The basic idea is to pull-back the structure of the unit ball of $L^a$ to construct a Banach manifold on the set of densities $\densities$. The actual construction is not straightforward in the infinite dimensional case, because the image set of the $\alpha$-embedding has empty interior in the Lebesgue space because of the non-negativity constraint $p^{1/a} \ge 0$. However, the approach works perfectly well in the cases of either a parametric model or a general model over a finite state space.

Variants of the $\alpha$-embedding has been studied by \cite{gibilisco|pistone:98}, \cite{burdet|combe|nencka:2001}, \cite{copas|eguchi:2005}. There is also literature covering the non-commutative (quantum) case, which is not treated here, e.g. \cite{gibilisco|isola:1999}, \cite{streater:2004orlicz}, \cite{jencova:2006}. 

The present paper is aimed to discuss what happens if we use $\Kln$ instead of the $\alpha$-embedding in the construction of the statistical manifold. The related issue of the algebro-statistical aspects of finite state space models is also discussed. See the books \cite{pistone|riccomagno|wynn:2001,pachter|sturmfels:2005,drton|sturmfels|sullivan:2009} on Algebraic Statistics and \cite{AGMS} on the relations between Algebraic Statistics and Information Geometry. 

Section 2 is an introduction to both the algebraic and the geometric features of the $\kappa$-deformed Gibbs model. Section 3 introduces a system of charts based on the $\kappa$-exponential; the resulting differentiable manifold is modeled on the space of centered and $1/\kappa$-integrable random variables. Section 4 discusses some properties of the tangent bundle; only a few basic results are presented. Section 5 is a short discussion on how the present construction relates to our formalism of exponential statistical manifolds and to other Authors' work. A full account the totality of the existing literature on Information Geometry is outside the scope of the present paper and will be published elsewhere.   
\section{$\kappa$-Deformed Gibbs model}
On a finite state space $\Omega$, equipped with the energy function $U \colon \Omega \rightarrow \reals_\ge$, we want to discuss the $\kappa$-deformation of the standard Gibbs model. There are two options, as there are two different presentation of the normalizing constant.

The first option would be to consider the statistical model
\begin{align}
  p(x;\theta) &= \frac{\kexpof{\theta U(x)}}{Z(\theta)} \label{eq:var-kgibbs} 
\\ &= \kexpof{\theta U(x) \kplus \klnof{\frac1{Z(\theta)}}} \notag
\end{align}
or, with $\widetilde\psi_\kappa(\theta) = \Kln Z(\theta)$,
\begin{multline}
  \Kln p(x;\theta) = \\ \theta U(x) \sqrt{1+\kappa^2 (\widetilde\psi_\kappa(\theta))^2} - \widetilde\psi_\kappa(\theta) \sqrt{1 + \kappa^2 \theta^2 U(x)^2}
\end{multline}
To the best of our knowledge, this model is not considered in the literature, therefore we are not discussing it here.

The second option is to define the generalised model as
\begin{align}
  p(x;\theta) &= \kexpof{\theta U(x) - \psi_\kappa(\theta)} \label{eq:kgibbs}\\
              &= \kexpof{\theta U(x)} \ktimes \kexpof{- \psi_\kappa(\theta)} \notag,
  \end{align}
where $\psi_\kappa(\theta)$ is the unique solution of the equation
\begin{equation}\label{eq:gibbs-psi}
  \sum_{x\in \Omega} \kexpof{\theta U(x) - \psi_\kappa(\theta)} = 1.
\end{equation}
This is the model that has been considered in \cite{kaniadakis:2002PhRE}, on the basis of the classical maximum entropy argument. It is a particular case of more general models as discussed by \cite{grunvald|dawid:2004} and by \cite{naudts:2004JIPAM,naudts:2008} and . We focus here on the algebraic and geometrical features of the model, not on its derivation from general principles.

The one-parameter statistical models \eqref{eq:var-kgibbs} and \eqref{eq:kgibbs} are different unless $\kappa=0$. This fact marks an important difference between the theory of ordinary exponential models and $\kappa$-deformed exponential models. For example, the derivative with respect to $\theta$ of the left hand side of \eqref{eq:gibbs-psi} is
\begin{equation}
  \sum_{x\in \Omega} \frac{U(x) - \psi_\kappa'(\theta)}{\sqrt{1+\kappa^2 \left(\theta U(x) - \psi_\kappa(\theta)\right)^2}} \kexpof{\theta U(x) - \psi_\kappa(\theta)}, 
\end{equation}
therefore
\begin{equation}
  \expectat \theta {\frac{U - \psi'(\theta)}{\sqrt{1+\kappa^2 \left(\theta U - \psi_\kappa(\theta)\right)^2}}} = 0,
\end{equation}
where $\expectat \theta V = \sum_x V(x) p(x;\theta)$. If $\kappa=0$, we have the usual formula $\psi'(\theta) = \expectat \theta U$.

From the geometrical point of view, the second approach has the advantage that the model for the $\Kln$-probabilities is linear. Here, by geometry we mean differential geometry of statistical models, i.e. the construction of an atlas of charts representing the subset of probability densities unto open sets of a Banach space, see \cite{lang:1995}. Before moving into this approach, we discuss the duality of the Gibbs model and its algebraic features. 

Let $V = \spanof{1,U}$ and $V^\perp$ the orthogonal space, i.e. $v \in V^\perp$ if, and only if,
$\sum_{x} v(x) = 0$ and $\sum_{x} v(x) U(x) = 0$. It follows from Equation \eqref{eq:kgibbs} that 
\begin{equation}\label{eq:kgibbsorth}
    \sum_{x \in {\Omega}} v(x) \klnof{p(x;\theta)} = 0, \qquad v \in V^\perp
\end{equation}
Vice versa, if a strictly positive probability density function $p$ is such that $\Kln p$ is orthogonal to $V^\perp$, then $p$ belongs to the $\kappa$-Gibbs model for some $\theta$.

For each $v \in V^\perp$, we can take its positive part $v^+$ and its negative part $v^-$, so that $v = v^+ - v^-$ and $v^+ v^- = 0$. Equation \eqref{eq:kgibbsorth} can be rewritten as 
\begin{equation}
\label{eq:kgibbsbin}
 \sum_{x\colon v(x) > 0} v^+(x) \klnof{p(x)} = \sum_{x\colon v(x) < 0} v^-(x) \klnof{p(x)}
\end{equation}
The interpretation of \eqref{eq:kgibbsbin} is the following. As $\sum_x v(x)=0$, we have
\begin{equation}
  \sum_{x \in {\Omega}} v^+(x) =   \sum_{x \in {\Omega}} v^-(x) = \lambda. 
\end{equation}
It follows that $r_1 = v^+/\lambda$, and $r_2 = v^-/\lambda$ are probability densities (states) with disjoint support, so that \eqref{eq:kgibbsbin} can be restated by saying that a positive density $p$ belongs to the $\kappa$-Gibbs model if, and only if,
\begin{equation}
  \expectat {r_1}{\klnof{p}} = \expectat {r_2}{\klnof{p}}
\end{equation}
for each couple of densities $r_1$, $r_2$ such that $r_1r_2=0$ and $\expectat {r_1} U = \expectat {r_2} U$, where $\expectat r V$ denotes the mean value of $V$ with respect to $r$, i.e. $\sum_x V(x) r(x)$.

If $v \in V^\perp$ happens to be integer valued, using the $\kappa$-algebra and the notation
\begin{equation}
\stackrel{\text{$n$ times}}{\overbrace{x \ktimes \cdots \ktimes x}} = x^{\ktimes n},
\end{equation}
we can write \eqref{eq:kgibbsbin} as
\begin{equation} \label{eq:kgibbs-zero}
 \stackrel[\kappa]{}{\bigotimes}_{x\colon v(x)>0} p(x)^{\ktimes v^+(x)} = {\stackrel[\kappa]{}{\bigotimes}_{x\colon v(x)<0} p(x)^{\ktimes v^-(x)}},
\end{equation}

It should be noted that \eqref{eq:kgibbs-zero} is continuous function of $p(x)$, $x \in \Omega$ and it does not require the strict positivity of each $p(x)$, $x \in \Omega$. Therefore, the same equation is satisfied by all limits (if any) of the model \eqref{eq:kgibbs}. Moreover, \eqref{eq:kgibbs-zero} is a \emph{$\kappa$-polynomial invariant} for the $\kappa$-Gibbs model, cf. in \cite{pistone:2008} the discussion of the case $\kappa=0$. Vice versa, each positive probability satisfying \eqref{eq:kgibbs-zero} belongs to the $\kappa$-Gibbs model.

The set of all $\kappa$-polynomial equations \eqref{eq:kgibbs-zero} is not a finite set, because each equation depends on the choice of an integer valued vector $v$ in the orthogonal space $V^\perp$. Accurate discussion of this issue requires tools from Commutative Algebra. If the energy function $U$ takes its values on a lattice $0, \Delta, 2\Delta, \dots$, $\Delta > 0$, we can choose integer valued random variables $v_1,\dots,v_{N-2}$ to be a linear basis of the orthogonal space $V^\perp$. In such a case, we have a finite system of binomial equations
\begin{equation}\label{eq:kgibbs-basis}
   \stackrel{\kappa}{\bigotimes}_{x\colon v_j(x)>0} p(x)^{\ktimes v^+_j(x)} =   \stackrel{\kappa}{\bigotimes}_{x\colon v(x)<0} p(x)^{\ktimes v^-_j(x)},
\end{equation}
$j = 1,\dots,N-2$, which is equivalent to the original model \eqref{eq:kgibbs}. However, it is not generally true that the finite system \eqref{eq:kgibbs-basis} is equivalent to the infinite system \eqref{eq:kgibbs-zero}.

In the classical case $\kappa=0$, the polynomial invariants of the Gibbs model form a polynomial ideal $I$ in the polynomial ring $\rationals[p(x)\colon x \in \Omega]$. The ideal I admits, because of the Hilbert Theorem, a finite generating set. The discussion of various canonical forms of such a generating set is one of the issues of Algebraic Statistics.

 \paragraph{Example} We specialize our discussion to the toy example which is discussed in \cite{pistone:2008} for $\kappa=0$. Let ${\Omega} = \set{1,2,3,4,5}$, $U(1) = U(2) = 0$, $U(3) = 1$, $U(4)= U(5) = 2$. The following display shows a set of integer valued $v_j$, $j=1,2,3$ of the orthogonal space $V^\perp$.
\begin{equation}
\bordermatrix[{[]}]{%
& 1 & U  & v_1 & v_2 & v_3 \cr
1& 1 & 0 & 1 & 0 & 1 & \cr
2& 1 & 0 & -1 & 0 & 1 \cr
3& 1 & 1 & 0 & 0 &  -4\cr
4& 1 & 2 & 0 & 1 & 1 \cr
5& 1 & 2 & 0 & -1 & 1 \cr
}%
\end{equation}
Equation \eqref{eq:kgibbs-basis} becomes:
\begin{equation} \label{eq:1}
  \left\{
    \begin{aligned}
      &p(1) = p(2) \\
      &p(4) = p(5) \\
      &p(1)\ktimes p(2)\ktimes p(4)\ktimes p(5) = p(3)^{\ktimes 4}
    \end{aligned}
\right.
\end{equation}
A strictly positive probability density belongs to the $\kappa$-Gibbs model \eqref{eq:kgibbs} if and only if it satisfies \eqref{eq:1}. The set of all polynomial invariants of the Gibbs model is a polynomial ideal and the question of finding a set of generators is intricate. A non strictly positive density that is a solution of \eqref{eq:1} is either $p(1)=p(2)=p(3)=0$, $p(4)=p(5)=1/2$, or $p(1)=p(2)=1/2$, $p(3)=p(4)=p(5)=0$. These two solutions are the uniform distributions on the sets of values that respectively maximize or minimize the energy function.

From \eqref{eq:15} we can derive another form of the last equation in the system \eqref{eq:1}:
\begin{multline}
 \left(p^\kappa(1) - \frac1{p^\kappa(1)}\right)+\left(p^\kappa(2) - \frac1{p^\kappa(2)}\right)+ \\ \left(p^\kappa(4) - \frac1{p^\kappa(4)}\right)+\left(p^\kappa(5) - \frac1{p^\kappa(5)}\right) = \notag \\ 4\left(p^\kappa(3) - \frac1{p^\kappa(3)}\right), \notag
\end{multline}
which is another algebraic form of the $\kappa$-Gibbs model. 

Jet another algebraic presentation is available. In the model \eqref{eq:kgibbs}, we introduce the new parameters
\begin{align}
  \zeta_0 &= \kexpof{- \psi_\kappa(\theta)}, \\
  \zeta_1 &= \kexpof{\theta},
\end{align}
so that
\begin{align}
  p(x;\theta) &= \kexpof{\theta U(x)} \ktimes \kexpof{- \psi_\kappa(\theta)}, \\
              &= \zeta_0 \ktimes \zeta_1^{\ktimes U(x)}.
\end{align}
The probabilities are $\kappa$-monomials in the parameters $\zeta_0, \zeta_1$, e.g.:
\begin{equation}\label{eq:toric-gibbs}
  \left\{
      \begin{aligned}
        p(1) &= p(2) = \zeta_0 \\ p(3) &= \zeta_0 \ktimes \zeta_1 \\ p(4) &= p(5) = \zeta_0 \ktimes \zeta_1^{\ktimes 2}
      \end{aligned}
\right.
\end{equation}

In algebraic terms, such a model is called a \index{toric model}\emph{toric model}. It is interesting to note that in \eqref{eq:toric-gibbs} the parameter $\zeta_0$ is required to be strictly positive, while the parameter $\zeta_1$ could be zero, giving rise the uniform distribution on $\set{1,2} = \set{x \colon U(x)=0}$. The other limit solution is not obtained by \eqref{eq:toric-gibbs}. The algebraic elimination of the indeterminates $\zeta_0, \zeta_1$ in \eqref{eq:toric-gibbs} will produce back polynomial invariants. For example, from $(\zeta_o\ktimes\zeta_1)^{\ktimes 2} = (\zeta_0)\ktimes(\zeta_0\ktimes\zeta_1^{\ktimes 2})$, we get $p(3)^{\ktimes 2} = p(2)\ktimes p(5)$. 

The monomial parametric representation of the type \eqref{eq:toric-gibbs}, together with their ability to parametrize the limit cases is discussed in \cite{pistone:2008}. The notion of toric model and the related discussion could be applied to the model \eqref{eq:var-kgibbs} which we are not discussing here.
\section{Charts}
\label{sec:chartintro}
We use the coordinate-free formalism of differential geometry of \cite{lang:1995}. Let us fix a $\kappa \in ]0,1[$. In order to construct an atlas of charts, i.e. mappings from the set of strictly positive probability densities $\pdensities$ to some model vector space, we consider each chart as associated to a density $p \in \pdensities$. Such a $p$ is used as a reference for every other density $q$ of a suitable subset of $\pdensities$, via the statistical notion of likelihood $q/p$. The theory is not restricted to a finite sample space. However, we do not discuss here the technicalities involved in the non finite case, but we use a general notation. 

We first define a $\kappa$-divergence. If $q, p \in \pdensities$ satisfy the condition
\begin{equation} \label{eq:100}
  \left(\frac qp\right)^\kappa, \left(\frac pq\right)^\kappa \in L^1(p),
\end{equation}
then the $\kappa$-divergence is defined to be
\begin{multline} \label{eq:divergence}
  D_\kappa(p \Vert q) = \expectat p {\klnof{\frac pq}} \\ = \frac1{2\kappa} \expectat p {\left(\frac pq\right)^\kappa - \left( \frac qp \right)^\kappa}.
\end{multline}
The first of the two conditions in \eqref{eq:100} is always satisfied because
\begin{equation}
  \label{eq:2}
  \expectat p {\left(\frac qp\right)^\kappa} \leq \left(\expectat p {\frac qp}\right)^\kappa = 1. 
\end{equation}
The second condition is non-trivial unless the state space is finite, because
we are assuming
\begin{equation}
  \label{eq:3}
\expectat p {\left(\frac pq\right)^\kappa} = \int \frac {p^{\kappa+1}}{q^\kappa} \,d\mu < +\infty.
\end{equation}
When such a condition is not satisfied the value of the expectation in \eqref{eq:divergence} is $+\infty$. We are not interested in this case here.

 The strict convexity of $-\Kln$ implies 
\begin{multline}
  D_\kappa(p \Vert q) = \expectat p {-\klnof{\frac qp}} \geq \\ -\klnof{\expectat p {\frac qp}} = \klnof 1 = 0.
\end{multline}
with equality if, and only if $q=p$.

The manifold we want to define has to be modelled on the Lebesgue space of centered $(1/\kappa)$-$p$-integrable random variables $L^{1/\kappa}_0(p)$, i.e. $v \in L^{1/\kappa}_0(p)$ if, and only if, $\expectat p {|v|^{1/\kappa}} = \int |v|^{1/\kappa}p\,d\mu < +\infty$ and $\expectat p v = \int vp\,d\mu=0$. At each $p$ there is a different model space, so that an isomorphism between them has to be provided. The simplest isometric identification between $L^{1/\kappa}(p_1)$ and $L^{1/\kappa}(p_2)$ is
\begin{equation} \label{eq:isometry}
  L^{1/\kappa}(p_1) \ni u \mapsto \left(\frac{p_1}{p_2}\right)^\kappa u \in L^{1/\kappa}(p_2)
\end{equation}
In fact,
\begin{equation}
  \label{eq:21}
  \int \left|\left(\frac{p_1}{p_2}\right)^\kappa u\right|^{1/\kappa} p_2\,d\mu = \int |u|^{1/\kappa} p_1\,d\mu. 
\end{equation}

We use a variation of the formalism used in the exponential case \cite{pistone|sempi:95}. The new manifold will be called $\kappa$-\emph{statistical manifold}. We define the subset $\mathcal E_p$ of $\pdensities$ by
\begin{align}
\mathcal E_{p} &= \setof{q \in \pdensities}{\left(\frac qp\right)^\kappa, \left( \frac pq\right)^\kappa \in L^{1/\kappa}(p)} \notag \\
              &= \setof{q \in \pdensities}{\frac qp, \frac pq \in L^1(p)} \notag \\ &= \setof{q \in \pdensities}{\frac pq \in L^1(p)} \label{eq:4}
\end{align}
The assumption in \eqref{eq:4} is stronger than that in \eqref{eq:3}, because here we want $\expectat p {p/q} = \int p^2/q \, d\mu < +\infty$, and $\expectat p {p^\kappa/q^\kappa} \le (\expectat p {p/q})^\kappa$. Indeed, we are assuming more than the mere existence of $D_\kappa(p \Vert q)$, unless the state space is finite. It should be noticed that the set $\mathcal E_p$ does not depend on $\kappa \in ]0,1[$. The condition $\int p^2/q \, d\mu < +\infty$ can be interpreted with the aid of the isometry \eqref{eq:isometry}. As
\begin{equation}
  \left(\frac pq\right)^\kappa \in L^{\frac1\kappa}(p)
\end{equation}
 by definition, we can map it to $L^{\frac1\kappa}(q)$ by multiplying it by itself, to get
 \begin{equation}
   \left(\frac pq\right)^{2\kappa} \in L^{\frac1\kappa}(q),\quad \text{or} \quad \left(\frac pq\right)^2 \in L^1(q).
 \end{equation}
Each of the sets $\mathcal E_p$, $p\in \pdensities$, is going to be the domain of a chart, and this will define an atlas of charts because of the covering $\pdensities = \cup_p \mathcal E_p$; a connected component of the manifold will be the union of overlapping $\mathcal E$'s.

If $q \in \mathcal E_{p}$, then $q$ is almost surely positive and we can write it in the form $q = \kexpof v p$, where
\begin{equation}
  v = \klnof{\frac qp} = \frac{\left(\frac qp\right)^{\kappa} - \left(\frac pq\right)^{\kappa}}{2\kappa} \in L^{1/\kappa}(p)
\end{equation}
The expected value at $p$ of $v = \klnof{\frac qp}$ is
\begin{equation}
\expectat p{\klnof{\frac qp}} = - D_{\kappa}(p \Vert q),
\end{equation}
 so that we can write every $q \in \mathcal E_{p}$ as
\begin{equation} \label{eq:7}
  q = \kexpof{u - D_\kappa(p \Vert q)} p,
\end{equation}
where $u$ is a uniquely defined element of the set of centered $(1/\kappa)$-$p$-integrable random variables $L^{1/\kappa}_0(p)$, namely
\begin{align}
  \label{eq:011}
  u &= \klnof{\frac qp} - \expectat p {\klnof{\frac qp}} \notag \\ &= \klnof{\frac qp} + D_\kappa(p \Vert q).
\end{align}

Conversely, given any $u \in L^{1/\kappa}_0(p)$, the real function $\psi \mapsto \expectat p {\kexpof{u - \psi}}$ is continuous and strictly decreasing from $+\infty$ to 0, therefore there exists a unique value of $\psi$, say $\psi_{\kappa,p}(u)$, such that
\begin{equation} \label{eq:8}
  \expectat p {\kexpof{u - \psi_{\kappa,p}(u)}} = 1,
\end{equation}
so that
\begin{equation} \label{eq:6}
  q = \kexpof{u - \psi_{\kappa,p}(u)} p \in \mathcal E_{p} \subset \pdensities.
\end{equation}

The 1-to-1 mapping
\begin{equation}
  \label{eq:5}
  \mathcal E_p \ni q \leftrightarrow u \in L^{1/\kappa}_0(p)
\end{equation}
is our chart. By comparing \eqref{eq:6} with \eqref{eq:7} we obtain $\psi_{\kappa,p}(u)=D_\kappa(p \Vert q)$, where $u$ is the image of $q$ in the chart at $p$. I.e. the $p$-chart representation of the functional $q \mapsto D_\kappa(p \Vert q)$ is $u \mapsto \psi_{\kappa,p}(u)$. A second example is the representation of $q \mapsto D_\kappa(q\Vert p)$:
\begin{equation}
  \label{eq:other-divergence}
  D_\kappa(q\Vert p) = \expectat q {\klnof{\frac qp}} = \expectat q u - \psi_{\kappa,p}(u).
\end{equation}

The functional $\psi_{\kappa,p}$ is of key importance in the case $\kappa=0$, where it is the cumulant functional of the random variable $u$:
\begin{equation}
  \psi_p(u) = \ln \expectat p {\euler^u}
\end{equation}

The proof of differentiability properties, in the non-finite case with $\kappa=0$ is not trivial, see e.g. \cite{cena|pistone:2007}, because the domain of  $\psi_p$ functional has a non-trivial description in that case. For $\kappa \ne 0$, the domain of $\psi_{\kappa,t}$ is the full space $L^{1/\kappa}_0(p)$, and we can compute the directional derivatives
\begin{equation}
  \label{eq:9}
  D \psi_{\kappa,p}(u) v = \left.\frac{d}{dt} \psi_{\kappa,p}(u + tv)\right|_{t=0}
\end{equation}
and
\begin{equation}
  \label{eq:800}
  D^2 \psi_{\kappa,p}(u) v w = \left.\frac{d^2}{ds dt} \psi_{\kappa,p}(u + sv + tw)\right|_{s=t=0}.
\end{equation}

For $\kappa=0$, the Fr\'echet derivatives of $\psi_p$ are
\begin{align}
  D \psi_p(u) v &= \expectat q v, \label{eq:Dpsi0}\\
  D^2 \psi_p(u) v w &= \covat q v w \label{eq:D2psi0},
\end{align}
where $v$ and $w$ are the directions of derivation and $q = \euler^{u - \psi_p(u)}  p$. 

For $\kappa \ne 0$ the computation of the directional derivative of \eqref{eq:8} gives
\begin{multline}
  \label{eq:90}
  \expectat p {\Kexp'\left(u - \psi_{\kappa,p}(u)\right)(v -D \psi_{\kappa,p}(u) v)} = \\ \expectat q {\frac{v -D \psi_{\kappa,p}(u) v}{\sqrt{1+\kappa^2(u - \psi_{\kappa,p}(u))^2}}} = 0.
\end{multline}
It follows from \eqref{eq:90} with $u=0$ that $D \psi_{\kappa,p}(0) = 0$. 

Otherwise, let $q|p$ denote the density proportional to 
\begin{multline}
  \label{eq:22}
  \Kexp'(u - \psi_{\kappa,p}(u)) p = \\ \frac{q}{\sqrt{1+\kappa^2(u - \psi_{\kappa,p}(u))^2}} \\ =   \frac{q}{\sqrt{1+\kappa^2\left(\klnof{\frac qp}\right)^2}},
\end{multline}
see \cite{naudts:2008}, where such a density is called escort probability. The explicit expression for the derivative is
\begin{equation}
  \label{eq:23}
  D \psi_{\kappa,p}(u) v = \expectat {q|p} v,
\end{equation}
which is the same as \eqref{eq:Dpsi0}, but the expectation is computed with respect of the escort density $q|p$. Later we will give a geometric interpretation of $q|p$.
 
The second derivative of $u \mapsto \kexpof{u - \psi_{\kappa,p}(u)}$ in the directions $v$ and $w$ is the first derivative in the direction $w$ of $u \mapsto \Kexp'(u - \psi_{\kappa,p}(u))(v - D\psi_{\kappa,p}(u) v)$, therefore it is equal to
\begin{multline}
  \label{eq:30}
\Kexp''(u - \psi_{\kappa,p}(u))(v-D\psi_{\kappa,p}(u)v)(w-D\psi_{\kappa,p}(u)w) \\ 
-\Kexp'(u - \psi_{\kappa,p}(u)) D^2\psi_{\kappa,p}(u) vw.
\end{multline}
The random variable in \eqref{eq:30} has zero $p$-expectation, so that 
\begin{multline}
  D^2\psi_{\kappa,p}(u) vw = \\ \frac{\expectat p {\Kexp''(u - \psi_{\kappa,p}(u))(v-D\psi_{\kappa,p}(u)v)(w-D\psi_{\kappa,p}(u)w)}}{\expectat p {\Kexp'(u - \psi_{\kappa,p}(u))}}.
\end{multline}
If $w=v\ne0$, then $D^2\psi_{\kappa,p}(u)vv > 0$, therefore the functional $\psi_{\kappa,p}$ is strictly convex. For $u=0$ we obtain
\begin{equation}
 D^2 \psi_{\kappa,p}(0) vw = \covat p u v.
\end{equation}
We do not have a similar interpretation for $u \ne 0$, but see the discussion of parallel transport below.
\section{The $\kappa$-statistical manifold and its tangent bundle}
\label{sec:kstat}
Assume now we want to change of chart, that is we want to change the reference density from $p$ to $\bar p$ to represent a $q$ that belongs both to $\mathcal E_{p}$ and to $\mathcal E_{\bar p}$. From now on, we skip the discussion of the non-finite case; such a discussion will be published elsewhere. In the finite state space case, the integrability conditions of \eqref{eq:4} are always satisfied, so that all the chart's domains $\mathcal E_p$ are equal to $\pdensities$. The application of \eqref{eq:011} and \eqref{eq:6} to the change
\begin{equation}
  \label{eq:300}
  L^{1/\kappa}_0(p) \ni u \mapsto q \mapsto \bar u \in L^{1/\kappa}_0(\bar p)
\end{equation}
gives
\begin{align}
  \bar u &= \klnof{\frac {q}{\bar p}} - \expectat {\bar p}{\klnof{\frac {q}{\bar p}}} \notag \\
      &= \klnof{\kexpof{u - \psi_{\kappa,p}(u)}\frac {p}{\bar p}} \notag \\ 
      &\qquad - \expectat{\bar p}{\klnof{\kexpof{u - \psi_{\kappa,p}(u)}\frac {p}{\bar p}}} \notag \\
      &= \left(u - \psi_{\kappa,p}(u)\right) \kplus \klnof{\frac {p}{\bar p}} \notag \\
      &\qquad - \expectat{\bar p}{\left(u - \psi_{\kappa,p}(u)\right) \kplus \klnof{\frac {p}{\bar p}}}
    \end{align}
For $\kappa=0$ the change of chart is an affine function:
\begin{equation}
  \label{eq:change0}
\bar u =  u + \lnof{\frac {p}{\bar p}} - \expectat{\bar p}{u + \lnof{\frac {p}{\bar p}}}
\end{equation}
with linear part
\begin{equation}
  \label{eq:13}
  v \mapsto v - \expectat {\bar p} v.
\end{equation}

For $\kappa \neq 0$, the derivative in the direction $v$ of the change of chart at $u$ is obtained from the derivation formula \eqref{eq:11}. It has the form $A - \expectat {\bar p} A$, with
\begin{multline}
  \label{eq:derivchange}
 A = \\ \left(\sqrt{1+\kappa^2 \left(\klnof{\frac {p}{\bar p}}\right)^2} - \frac{\kappa^2 \klnof{\frac {p}{\bar p}} \left(u - \psi_{\kappa,p}(u)\right)}{\sqrt{1+\kappa^2 \left(u - \psi_{\kappa,p}(u)\right)^2}}\right) \\ \times\left(v - D \psi_{\kappa,p}(u) v\right).
\end{multline}
      
We are then led to the study the tangent spaces of the $\kappa$-statistical manifold. Let $p_\theta$, $\theta \in ]-1,1[$, be a curve in $\mathcal E_{p}$,
    \begin{equation} \label{eq:curve}
      p_\theta = \kexpof{u_\theta -\psi_{\kappa,p_0}(u_\theta)} p_0.
    \end{equation}
In the chart at $p$ the velocity vector is given by 
    \begin{equation}
      \dot{u_\theta} \in L^{1/\kappa}_0(p).
    \end{equation}

We identify the tangent space at $p_0$ with the space of the random variables $\dot{u_0}$. In general, the tangent space at $p \in \pdensities$ is defied to be $T_p = L^{1/\kappa}_0(p)$. Derivation with respect to $\theta$ of \eqref{eq:curve} gives
  \begin{equation}
\label{eq:velocity}
    \frac{\dot{p_\theta}}{p_\theta} = (1 + \kappa^2(u_\theta -\psi_{\kappa,p}(u_\theta)^2)^{-1/2} (\dot{u_\theta} -D\psi_{\kappa,p}(u_\theta)(\dot{u_\theta})).
  \end{equation}
In particular, $\dot{p_0}/p_0 = \dot{u_0}$. To extend this to the general $\theta$, let us observe first that \eqref{eq:90}, with $p=p_0$, $u = u_\theta$, $v=\dot{u_\theta}$, shows that $\expectat {p_0}{\dot{p_\theta}/p_\theta} = 0$. The second factor in the right end side of \eqref{eq:velocity} is, because of \eqref{eq:23},  
\begin{equation}
  \dot{u_\theta} -D\psi_{\kappa,p}(u_\theta)(\dot{u_\theta}) = \dot{u_\theta} - \expectat {p_\theta|p_0}{\dot{u_\theta}}.
\end{equation}
The first factor is
\begin{multline}
  \left(1 + \kappa^2(u_\theta -\psi_{\kappa,p}(u_\theta)^2\right)^{-1/2} = \\ \left(1 + \kappa^2\left(\klnof{\frac{p_\theta}{p_0}}\right)^2\right)^{-1/2}
\end{multline}

Let us define a parallel transport $U^\kappa_{p,\bar p}$ mapping the tangent space at $p$, i.e. $T_{p} = L^{1/\kappa}_0(p)$, on the tangent space at $\bar p$, i.e. $T_{\bar p} = L^{1/\kappa}_0(\bar p)$ as follows. If $u \in L^{1/\kappa}_0(p)$, then the random variable 
\begin{equation}
  \label{eq:connection}
  \bar u = \frac{u - \expectat {\bar p|p}{u}}{\sqrt{1 + \kappa^2\left(\klnof{\frac {\bar p}{p}}\right)^2}}
\end{equation}
belongs to $T_{\bar p} = L^{1/\kappa}_0(\bar p)$, and we define $\bar u = U^\kappa_{p,\bar p}(u)$. In fact, $\expectat {\bar p} {\bar u} = 0$, and, moreover, 
\begin{equation}
  \left(1 + \kappa^2\left(\klnof{\frac {\bar p}{p}}\right)^2\right)^{-1/2} \le 2 \left(\frac {\bar p}{p}\right)^\kappa,
\end{equation}
so that the conclusion 
\begin{equation}
  \expectat {\bar p}{\vert U^\kappa_{p,\bar p}(u)\vert^{1/\kappa}} \le 2^{1/\kappa} \expectat {p}{\vert u\vert^{1/\kappa}}
\end{equation}
follows from the isometry \eqref{eq:isometry}.

The quantity $\dot{p_\theta}/p_\theta$ in \eqref{eq:velocity} has the following remarkable interpretation in terms of the parallel transport:
\begin{equation} \label{eq:parallel}
  \frac{\dot{p_\theta}}{p_\theta} = U^\kappa_{p_0,p_\theta} (\dot{u_\theta}),
\end{equation}
i.e. it is the velocity vector at $\theta$, transported to the tangent space at $p_\theta$.

In this setting, we can define $\kappa$-exponential models in a non parametric way as 
\begin{equation}\label{eq:kexpmodel}
 q = \kexpof{u - \psi_{\kappa,p}(u)} p, \quad u \in V,
\end{equation}
where $V$ is a linear sub-space of $L_0^{1/\kappa}(p)$. Each $v \in V$ is called a canonical variable of the $\kappa$-exponential model. The implicit representation of the exponential model \eqref{eq:kexpmodel} is
\begin{equation}
\expectat p {\klnof{\frac qp} v} = 0, \quad v \in V^\perp .
\end{equation}
where $V^\perp \subset L^{1-1/\kappa}_0(p)$ is the orthogonal of $V$. In Statistics, a $v \in V^\perp$ is called a constraint of the log-linear model. We could derive, as we did in the $\kappa$-Gibbs model example, for lattice-valued constraint variables $v$, the relevant polynomial-type equations based of the deformed product operation. 

In particular, a one dimensional exponential model is characterised by a one-dimensional space $V = \spanof{u}$, $u \in L^{1/\kappa}_0(p)$. Let us show that such a model satisfies a differential equation on the manifold. Given $u \in L^{1/\kappa}_0(p_0)$, for each $q \in \mathcal E_p$, we can define the mapping
\begin{equation}
  \label{eq:field}
  \mathcal E_p \ni q \mapsto U^\kappa_{p,q}(u) \in T_q,
\end{equation}
which is a vector field of the $\kappa$-manifold. The velocity of a curve $p_\theta=\kexpof{u_\theta - \psi_{\kappa,p_0}(u_\theta)} p_0$ is represented in the chart at $p_\theta$ by $\dot{p_\theta}/{p_\theta}$ because of \eqref{eq:parallel}, therefore we can consider the differential equation
\begin{equation}
  \label{eq:differential}
  \frac{\dot{p_\theta}}{p_\theta} =  U^\kappa_{p_0,p_\theta}(u),
\end{equation}
whose solution is
\begin{equation}
  p_\theta = \kexpof{\theta u - \psi_{\kappa,p}(\theta u)} p_0,
\end{equation}
cf. \eqref{eq:velocity}. 
Equation \eqref{eq:differential} implies that the one dimensional exponential model has tangent vectors which are transported into each other by the parallel transport $U^\kappa$, i.e. the model is $U^\kappa$-auto-parallel.
 
The one dimensional exponential model is one of the simplest example of differential equation on the $\kappa$-statistical manifold. The treatment of evolution equations for densities is one of the main motivation for introducing a manifold structure on the set of densities, cf. \cite{otto:2001}, \cite{ay|erb:2005}, \cite{ohara|wada:2008}.
\section{Statistical manifolds}
In this final section, we go back to the general setting and briefly discuss how the $\kappa$-statistical manifolds we have defined relate with the previous construction of the exponential statistical manifold.
$(\Omega,\mathcal F, \mu)$ is a generic probability space, $\sdensities$ is the set of real random variables $f$ such that $\int f\,d\mu =1$, $\densities$ the convex set of probability densities, $\pdensities$ the convex set of strictly positive probability densities: 
\begin{equation}
  \pdensities \subset \densities \subset \sdensities
\end{equation}

In the classical case, i.e. $\kappa=0$, differentiable manifolds are defined on both $\pdensities$ and $\sdensities$. Such manifolds are both modeled on suitable Orlicz spaces, see \cite{rao|ren:2002}. Orlicz spaces are a generalization of Lebesgue spaces, where the norm is defined through a symmetric, null at zero, non-negative, convex, with more than linear growth, function called Young function. Let $\Phi$ be any Young function with growth equivalent to $\exp$, e.g. $\Phi(x) = \cosh(x)-1$, with convex conjugate $\Psi$, e.g. $\Psi(y) = (1+\absoluteval y) \logof{1+\absoluteval y} - \absoluteval y$. The relevant Orlicz spaces are denoted by $L^{\Phi}$ and $L^{\Psi}$, respectively. Both these Banach spaces appear naturally in Statistics. A random variable $u$ belongs to the space $L^\Phi$ if, an only if, its Laplace transform exists in a neighborhood of 0, that is, the one dimensional exponential model $p(\theta) \propto \euler^{\theta u}$ is defined for values of the parameter in a neighborhood of 0. A density function $f$ has finite entropy if, and only if, it belongs to the space $L^\Psi$. We denote by $L_0^{\Phi}$, $L_0^{\Psi}$ the sub-spaces of centered random variables. If the sample space is not finite, then the exponential Orlicz space is not separable and the closure $M^{\Phi}$ of the space of bounded functions is different from $L^{\Phi}$. There is a natural separating duality between $L^{\Phi}_0$ and $L^{\Psi}_0$, which is given by the bi-linear form
\begin{equation}
  (u,v) \mapsto \int uv\,d\mu.
\end{equation}
It can be proved that the cumulant generating functional $\psi_p(u) = \expectat p {\euler^u}$, $u \in \L^\Psi_0(p)$, $p \in \pdensities$ is positive, strictly convex, analytic. The interior of the proper domain of $\psi_p$,
\begin{equation}
  \mathcal S_p = \setof{u \in \L^\Psi_0(p)}{\psi_p(u) < +\infty}^\circ
\end{equation}
 defines the so called maximal exponential model
 \begin{equation}
   \mathcal E_{0,p} = \setof{\euler^{u-\psi_p(u)} p}{u \in \mathcal S_p}
 \end{equation}

In our case, i.e. $0 < \kappa  < 1$, we could consider
\begin{multline}
  \kcoshof x - 1 = \frac12(\kexpof{x}+\kexpof{-x}) - 1 \\
             = \frac12\left(\euler^{\int_0^x \frac{dt}{\sqrt{1+\kappa^2 t^2}}}+\euler^{\int_0^{-x} \frac{dt}{\sqrt{1+\kappa^2 t^2}}}\right) - 1 \\
             = \cosh\left(\int_0^x \frac{dt}{\sqrt{1+\kappa^2 t^2}}\right) - 1.
\end{multline}
We obtain a Young function equivalent to $|x|^{1/\kappa}$, so that the related Orlicz space is just $L^{1/\kappa}$ and the $\kappa$-statistical manifolds appears as a generalisation of the exponential construction.

If $q \in \mathcal E_{0,p}$, then $q = \expof{u - \psi_p(u)} . p$. If $-u \in \mathcal S_p$, then $q \in \mathcal E_p$, so that $q = \kexpof{v - \psi_{k,p}(v)}$ for a suitable $v \in L^{1/\kappa}_0(p)$. It follows that
\begin{equation}
  \expof{u - \psi_p(u)} = \kexpof{v - \psi_{k,p}(v)},
\end{equation}
therefore the e-coordinate $u$ and $\kappa$-coordinate $v$ are related by the equation
\begin{equation}
  u - \psi_p(u) = \int_0^{v - \psi_{k,p}(v)} \frac{dt}{\sqrt{1 + \kappa^2 t^2}}.
\end{equation}
\section{Conclusion}
We have presented a non parametric construction of the statistical manifold based on the use of the centered $\Kln$-likelihood as a functional coordinate. The basics of the formalism are discussed in Section 3 and 4. A notable difference from the standard case $\kappa=0$ is the absence of simple formul{\ae} for the functional $\psi_{\kappa,p}$ and its derivatives. It has been shown that the derivation of such quantities is related with a suitably defined parallel transport of the tangent bundle; for such a parallel transport the exponential models are auto-parallel. This, in turn, should lead to a theory of evolution equations on the $\kappa$-manifold. The algebraic features of the $\kappa$-exponential models for lattice contrasts have been described on an example. The relation of this construction to $\alpha$-geometries of \cite{amari|nagaoka:2000} is not discussed here nor is the relation with other geometric construction based on different divergences, e.g. Tsallis entropy \cite{tsallis:1988}. When comparing our formalism with similar construction of the statistical manifold based on other transformations derived from the power transform, it should be noticed that the $\kappa$-logarithm has a distinctive advantage of having range $\reals$ while retaining a simple algebraic character.
\begin{acknowledgement} 
This piece of research was started in occasion of the SigmaPhi 2008 conference, Kolympari, GR. The Author thanks the organizers of the symposium Fisher Information and Geometry for the invitation to present a paper there and to R. Trasarti-Battistoni for the conversations during the conference. A preliminary version of this paper has been discussed with G. Kaniadakis and was presented at MSRI Berkeley CA. The work has been supported by DIMAT, Politecnico di Torino, by MSRI Berkeley CA, and by SAMSI, Research Triangle Park NC.
\end{acknowledgement}
%
%\bibliography{/Users/gianni/Archive/bibs/tutto}% 

%
\end{document}